\documentclass{article}

\usepackage[backend=bibtex,maxbibnames=99]{biblatex}
\renewbibmacro{in:}{}
\usepackage{xcolor}
\usepackage[margin=1in]{geometry}
\usepackage{amsfonts}
\usepackage{amssymb}
\usepackage{amsmath}
\usepackage{bm}
\usepackage{fancyhdr}
\usepackage{tikz}
\usepackage{ytableau}
\usepackage[utf8]{inputenc}
\usepackage{enumitem}
\usepackage{amsthm}
\usepackage{thmtools,thm-restate}
\usepackage{multirow}

\newcommand{\C}{\mathbb{C}}

\let\P\relax
\newcommand{\P}{\mathbb{P}}

\DeclareMathOperator{\SSYT}{SSYT}
\DeclareMathOperator{\Sym}{Sym}
\DeclareMathOperator{\Hom}{Hom}
\DeclareMathOperator{\coker}{coker}
\DeclareMathOperator{\Ima}{Im}
\DeclareMathOperator{\rk}{rk}

\DeclareMathOperator{\Id}{Id}

\newtheorem{prop}{Proposition}
\newtheorem{theorem}{Theorem}
\newtheorem{cor}{Corollary}

\theoremstyle{remark}
\newtheorem*{remark}{Remark}

\theoremstyle{definition}
\newtheorem{example}{Example}

\begin{document}

\title{Border rank bounds for $GL(V)$-invariant tensors arising from matrices of constant rank}

\author{Derek Wu\thanks{Wu partially supported by NSF grant  AF-2203618}}

\date{}

\maketitle

\begin{abstract}
  We prove border rank bounds for a class of $GL(V)$-invariant tensors in $V^*\otimes U\otimes W$, where $U$ and $W$ are $GL(V)$-modules.
  These tensors correspond to spaces of matrices of constant rank.
  In particular we prove lower bounds for tensors in $\C^l\otimes\C^m\otimes\C^n$ that are not $1_A$-generic, where no nontrivial bounds were known, and also when $l,m\ll n$, where previously only bounds for unbalanced matrix multiplication tensors were known.
  We give the first explicit use of Young flattenings for tensors beyond Koszul to obtain border rank lower bounds, and determine the border rank of three tensors.
\end{abstract}

\section{Introduction}

Let $A$, $B$, and $C$ be complex vector spaces, of dimensions $l$, $m$, and $n$ respectively.
A tensor $T\in A\otimes B\otimes C$ has \emph{rank one} if $T=a\otimes b\otimes c$ for some $a\in A$, $b\in B$, and $c\in C$.
The \emph{rank} of $T$ is the smallest integer $r$ such that $T$ is the sum of $r$ rank one tensors; we write $R(T)=r$.
The \emph{border rank} $\underline{R}(T)$ is the smallest integer $r$ such that $T$ is the limit of a sum of $r$ rank one tensors.
We have $R(T)\geq \underline{R}(T)$, and this inequality can be strict.

Finding, or even bounding, the border rank of a tensor is an open problem, with applications in theoretical computer science \cite{15}, algebraic statistics \cite{14}, and tensor decomposition \cite{16}.
Several lower bound techniques have been developed in the last ten years, including Young flattenings \cite{3,2}, border substitution \cite{17,18}, and border apolarity \cite{12,1}.
This paper advances these techniques, and applies them to an interesting class of tensors.

A tensor is \emph{concise} if the maps $T_A\colon A^*\to B\otimes C$, $T_B$, and $T_C$ (defined similarly) are all injective.
To avoid trivialities, all tensors we consider in this paper are concise.
If a tensor $T$ is concise, then $\underline{R}(T)\geq\max\{l,m,n\}$.
If a concise tensor $T$ satisfies $R(T)=\max\{l,m,n\}$, we say $T$ is of \emph{minimal rank}.
Similarly if $T$ satisfies $\underline{R}(T)=\max\{l,m,n\}$, then $T$ is of \emph{minimal border rank}.
A tensor is \emph{$1_A$-generic} if there exists $\alpha\in A^*$ such that $T_A(\alpha)\in B\otimes C$ is of full rank $\min\{m,n\}$.
Tensors following the construction in Section 2.2 are not $1_A$-generic.

The set of rank at most one tensors in $A\otimes B\otimes C$ is the cone over the Segre variety.
The set of tensors of border rank at most $r$ is the cone over the $r$-th secant variety of the Segre variety, so in particular it is an algebraic variety defined by polynomials.
For a tensor $T$, one way of showing $\underline{R}(T)>r$ is by finding a polynomial that vanishes on the $r$-th secant variety, but is nonzero when evaluated at $T$.

In 1983, Strassen \cite{6} found equations for lower bounds up to $\frac{3}{2}m$ for tensors in $\C^k\otimes\C^m\otimes\C^m$ when $k\geq3$.
In 2013, Landsberg and Ottaviani \cite{3,2} used representation theory to find equations given by what they called \emph{Young flattenings}.
They showed Young flattenings give nontrivial equations for \emph{Waring border rank} of symmetric tensors in \cite[Lemma~4.2.1,Proposition~4.2.5]{3}.
In \cite{2}, they used \emph{Koszul flattenings}, a special case of Young flattenings, to obtain lower bounds of $2m-\log m$ for tensors in $\C^k\otimes\C^m\otimes\C^m$, with $k\leq m$.
In particular, they improved the border rank lower bound of the matrix multiplication tensor.

The tensors we study originate from \cite[Section~4]{4}, where Landsberg and Manivel construct $GL(V)$-invariant spaces of matrices of constant rank.
Let $V$ be a $k$-dimensional vector space over $\C$, and $\mu$ a partition with at most $k$ nonzero parts.
We may identify $\mu=(\mu_1,\ldots,\mu_k)$ with its \emph{Young diagram}, a collection of left and top aligned \emph{cells}, with $\mu_j$ cells in the $j$-th row.
Let $\nu$ be a partition obtained by adding a cell to $\mu$ such that we still have a Young diagram after adding the cell.
The partitions $\mu$ and $\nu$ give rise to $GL(V)$-modules $S_\mu V^*$ and $S_\nu V$, of dimension $m$ and $n$ respectively.
We obtain a tensor $T\in V^*\otimes S_\mu V^*\otimes S_\nu V\cong\C^k\otimes\C^m\otimes\C^n$ that corresponds to the $GL(V)$-equivariant map
\begin{gather*}
  \phi\colon V\to \Hom(S_\mu V,S_\nu V)\text{,}
\end{gather*}
a space of matrices of constant rank (see Section 2.2).
If $\nu$ differs from $\mu$ in the first or $k$-th row, the space of matrices is full rank, and otherwise it is of bounded rank.
For $v\in V$, denote the image of $v$ under $\phi$ by $\phi_v$.
The map is $\phi_v(u)=P(u\otimes v)$, where $P\colon S_\mu V\otimes V\to S_\nu V$ is the $GL(V)$-equivariant projection map.

Previously, the only known nontrivial bound for tensors $T\in A\otimes B\otimes C$, $l\ll m<n$ is that $T$ is not minimal border rank if $T_C(C^*)\subseteq A\otimes B$ does not contain a rank one element \cite[Exercise~15.14]{19}.
Additionally, previous results for border rank lower bounds for tensors in $A\otimes B\otimes C$ with $l,m\ll n$ were only for the matrix multiplication tensors $M_{\langle2nn\rangle}\in \C^{2n}\otimes\C^{2n}\otimes\C^{n^2}$ and $M_{\langle3nn\rangle}\in \C^{3n}\otimes\C^{3n}\otimes\C^{n^2}$ \cite{1}.
We give nontrivial border rank bounds for a class of tensors with larger imbalances in dimension than $M_{\langle2nn\rangle}$ and $M_{\langle3nn\rangle}$.

\subsection{Results and conjectures}

First, we consider the case where $\nu$ is obtained from $\mu$ by adding to the first or $k$-th row of $\mu$; we state the result in full generality. Unlike other tensors considered in this paper, the these tensors do not correspond to spaces of matrices of bounded rank.

For a reductive group $G$ and $G$-modules $A$ and $B$ with highest weight $\lambda$ and $\mu$ respectively, the \emph{Cartan product} $C$ is the $G$-module with highest weight $\lambda+\mu$.
The module $C$ occurs with multiplicity one in $A\otimes B$.

\begin{restatable}{prop}{Cprod}\label{Cprod}
  Let $A$ and $B$ be irreducible modules of a reductive group $G$, and $C$ the Cartan product of $A$ and $B$.
  Then, the unique $G$-invariant tensor $T\in A^*\otimes B^*\otimes C$ is of minimal rank $\dim C$.
\end{restatable}

We give one example of an explicit rank decomposition.

For the remaining results, we consider the case where $\nu$ is obtained from $\mu$ by adding a cell to neither the first nor $k$-th row.
Because of $GL(V)$-invariance, Young flattenings of these tensors are of particular interest: all our flattening maps are $GL(V)$-module maps.
Our effective use of Young flattenings beyond Koszul flattenings is the first such for border rank of tensors, and are specifically adapted to the tensors we consider.
In some cases where the ambient space is $\C^l\otimes\C^m\otimes\C^n$, $m\neq n$, general Young flattenings can give better bounds than those given by Koszul flattenings.

With our tensor as the map $T\colon S_\nu V^*\to V^*\otimes S_\mu V^*$, the \emph{$p$-th Koszul flattening} is the composition of $GL(V)$-equivariant maps
\begin{gather}\label{equi}\tag{$\ast$}
  T^{\wedge p}\colon \Lambda^pV^*\otimes S_\nu V^*
  \to \Lambda^pV^*\otimes V^*\otimes S_\mu V^*
  \to \Lambda^{p+1}V^*\otimes S_\mu V^*\text{.}
\end{gather}
This map gives a border rank lower bound of
\begin{gather*}
  \underline{R}(T)\geq\left\lceil\frac{\rk(T^{\wedge p})}{\binom{k-1}{p}}\right\rceil\text{,}
\end{gather*}
where $\rk(T^{\wedge p})$ is the rank of the linear map (\ref{equi}).
We state the following theorem, which says that all modules that can appear in the image, subject to Schur's lemma, do appear:

\begin{restatable}[{\cite[Lemma~2.1,Remark~2.4]{13}}]{theorem}{pKoszul}\label{pKoszul}
  Let $\mu$ and $\nu$ be two partitions, such that $\nu$ is obtained by adding a cell in row $j$ to $\mu$.
  Then, the image of the $p$-th Koszul flattening map
  \begin{gather*}
    T^{\wedge p}\colon
    \Lambda^{p}V^*\otimes S_\nu V^*
    \to \Lambda^{p+1} V^*\otimes S_\mu V^*
  \end{gather*}
  contains all submodules appearing in both the domain and codomain.
  These are given by the partitions $\pi$ with $\pi_j=\mu_j+1=\nu_j$ for $S_\pi V\in\Lambda^p V\otimes S_\nu V$.
\end{restatable}

The Koszul flattening maps are exactly the differential maps appearing in \emph{Pieri resolutions} over exterior algebras, which can be found in \cite{13}.

Introduce the \emph{Young flattening} of $T\in V^*\otimes S_\mu V^*\otimes S_\nu V$ with partitions $\alpha$ and $\tilde{\alpha}$, where $\tilde{\alpha}$ obtained by adding a cell to $\alpha$.
It is the composition of  $GL(V)$-equivariant maps
\begin{gather*}
  T_{\alpha,\tilde{\alpha}}
  \colon S_\alpha V^*\otimes S_\nu V^*
  \to S_\alpha V^*\otimes V^*\otimes S_\mu V^*
  \to S_{\tilde{\alpha}}V^*\otimes S_{\mu}V^*
\end{gather*}
(see Sections 2.3, 3.2 for an explicit description).
We note the $p$-th Koszul flattening is the special case $\alpha=(1^p)$ and $\tilde{\alpha}=(1^{p+1})$.
The map $T_{\alpha,\tilde{\alpha}}$ gives a border rank lower bound of
\begin{gather*}
  \left\lceil\frac{\rk(T_{\alpha,\tilde{\alpha}})}{r}\right\rceil
\end{gather*}
where $r$ is the rank of a generic matrix in the space
\begin{gather*}
  V^*\subset\Hom(S_{\alpha}V^*,S_{\tilde{\alpha}}V^*)\text{.}
\end{gather*}
We make the following observation: if we choose $\alpha=\mu$ and $\tilde{\alpha}=\nu$, then the resulting flattening map, which we denote $\Phi$, has isomorphic domain and codomain.
Thus, by $GL(V)$-invariance, $\Phi$ may be an isomorphism, and we conjecture this holds:

\begin{restatable}{conj}{Yflat}\label{Yflat}
  Let $V$ be a $k$-dimensional vector space, and $\mu$ and $\nu$ two partitions such that $\nu$ is obtained by adding a cell to $\mu$ not in the first or $k$-th row.
  Then the map
  \begin{gather*}
    \Phi\colon
    S_\mu V^*\otimes S_\nu V^*
    \to S_\mu V^*\otimes V^*\otimes S_\mu V^*
    \to S_\nu V^*\otimes S_\mu V^*
  \end{gather*}
  given by inclusion into the intermediate space, followed by a projection to the final space, is of full rank.
\end{restatable}

In particular, proving this conjecture would show that no tensor obtained via the construction in Section 2.2 is of minimal border rank.
Since the nonzero matrices $V\subset\Hom(S_\mu V,S_\nu V)$ are of bounded constant rank $r$, we have
\begin{gather*}
  r<\min\{m,n\}\text{.}
\end{gather*}
If the Young flattening is full rank, then
\begin{gather*}
  \underline{R}(T)\geq\left\lceil\frac{mn}{r}\right\rceil>\max\{m,n\}\text{.}
\end{gather*}
Because of the bounded rank condition, we require that $\nu$ cannot be obtained from $\mu$ by adding to the first or $k$-th row.
If $\nu$ is obtained from $\mu$ by adding to the first or $k$-th row, we use Proposition 1.

We have partial results for this conjecture:
\begin{restatable}{prop}{isotypic}\label{isotypic}
  If an isotypic component with highest weight $\pi$ appears with the same multiplicity $M$ in $S_\mu V^*\otimes S_\nu V^*$ and $S_\mu V^*\otimes V^*\otimes S_\mu V^*$, then the isotypic component of $\pi$ appears in the image of the flattening map $T'\colon S_\mu V^*\otimes S_\nu V^*\to S_\nu V^*\otimes S_\mu V^*$ with multiplicity $M$.
\end{restatable}

For small partitions we verify the conjecture computationally.

\begin{restatable}{prop}{smalltensors}\label{smalltensors}
  If $\nu$ fits inside a $4\times 4$ square or $5\times 2$ rectangle, then the Young flattening map
  \begin{gather*}
    T'\colon
    S_\mu V^*\otimes S_\nu V^*
    \to S_\nu V^*\otimes S_\mu V^*
  \end{gather*}
  is of full rank.
\end{restatable}

Thus any $GL(V)$-invariant tensor $T\in V^*\otimes S_\mu V^*\otimes S_\nu V$ with $\nu$ fitting inside a $4\times 4$ square or $5\times 2$ rectangle is not of minimal border rank for $k>\ell(\nu)$.
Once $T'$ is full rank for $\dim V=k$, $k$ large enough such that all $GL(V)$-modules in $S_\nu V\otimes S_\mu V$ are nonzero, then $T'$ is full rank for all $\dim V\geq k$.

For most spaces in this paper there are tensors of far larger border rank than the lower bounds we obtain.
This is typical of the state of the art, and it is difficult to determine when these bounds are sharp.
However, for the $GL_3$-invariant tensor $T\in (\C^3)^*\otimes S_{(2,1)}(\C^3)^*\otimes S_{(2,2)}\C^3$, the dimension of its ambient space is small enough for the Young flattening to determine is border rank.

\begin{restatable}{theorem}{tse}\label{tse}
  The $GL_3$-invariant tensor $T\in (\C^3)^*\otimes S_{(2,1)}(\C^3)^*\otimes S_{(2,2)}\C^3\cong \C^3\otimes\C^8\otimes\C^6$ has maximal border rank 10.
\end{restatable}

The first Koszul flattening only gives a border rank lower bound of nine.

We give explicit bounds for $\dim V=3$, and $\mu,\nu$ 2-row partitions in a $5\times2$ rectangle.

\begin{restatable}{theorem}{smallab}\label{smallab}
  Let $\dim V=3$, $\mu=(a+b+1,a)$, and $\nu=(a+b+1,a+1)$ for $a,b\geq0$.
  Border rank lower bounds for the tensors $T_{a,b}\in V^*\otimes S_{\mu}V^*\otimes S_{\nu}V\cong \C^3\otimes\C^m\otimes\C^n$ are given in the following table:
  \begin{center}
  \begin{tabular}{c|c c c c c}
  $\underline{R}(T_{a,b})\geq$ & $b=0$ & $1$ & $2$ & $3$ & $4$ \\
  \hline
  $a=0$ & $5$ & $10$ & $17$ & $26$ & $37$\\
  $1$ & $10$ & $19$ & $31$ & $46$\\
  $2$ & $17$ & $31$ & $49$\\
  $3$ & $26$ & $46$\\
  $4$ & $37$
\end{tabular}\text{.}
  \end{center}
\end{restatable}

For comparison, we provide a table for $m=\dim S_\mu V$ and $n=\dim S_\nu V$:

\begin{center}
\begin{tabular}{c|c c c c c}
  $(m,n)$ & $b=0$ & $1$ & $2$ & $3$ & $4$ \\
  \hline
  $a=0$ & $(3,3)$ & $(6,8)$ & $(10,15)$ & $(15,24)$ & $(21,35)$\\
  $1$ & $(8,6)$ & $(15,15)$ & $(24,27)$ & $(35,42)$\\
  $2$ & $(15,10)$ & $(27,24)$ & $(42,42)$\\
  $3$ & $(24,15)$ & $(42,35)$\\
  $4$ & $(35,21)$
\end{tabular}\text{,}
\end{center}

and a table for the border rank lower bound given by a Koszul flattening:

\begin{center}
\begin{tabular}{c|c c c c c}
  $\underline{R}(T_{a,b})\geq$ & $b=0$ & $1$ & $2$ & $3$ & $4$ \\
  \hline
  $a=0$ & $5$ & $9$ & $15$ & $23$ & $32$\\
  $1$ & $9$ & $18$ & $29$ & $42$\\
  $2$ & $15$ & $29$ & $46$\\
  $3$ & $23$ & $42$\\
  $4$ & $32$
\end{tabular}\text{.}
\end{center}

Next, let $T_k\in V^*\otimes V^*\otimes \Lambda^2 V\cong\C^k\otimes\C^k\otimes\C^{\binom{k}{2}}$ denote the tensor corresponding to the space of skew-symmetric matrices.
The first Koszul flattening gives the best border rank lower bounds, no matter the dimension of $V$.
Because of the imbalance of dimension, we do not obtain better bounds by restricting to a subspace; see Section 5.1.
We also use explicit border rank decompositions to obtain upper bounds.

\begin{theorem}
  Let $T_k\in V^*\otimes V^*\otimes \Lambda^2 V\cong\C^k\otimes\C^k\otimes\C^{\binom{k}{2}}$ be the unique $GL(V)$-invariant tensor.
  For $k\geq4$,
  \begin{gather*}
    \left\lceil\frac{k^2}{2}\right\rceil
    = \binom{k}{2}+\left\lceil\frac{k}{2}\right\rceil
    \leq\underline{R}(T_k)
    \leq \binom{k}{2}+k-2
    =\binom{k+1}{2}-2\text{.}
  \end{gather*}
\end{theorem}

It was previously known that $\underline{R}(T_3)=5$ (see \cite[Proposition 3.2]{7}), and these bounds give $\underline{R}(T_4)=8$ and $\underline{R}(T_5)=13$.

Finally, we apply techniques of \emph{border apolarity} \cite{12} following \cite{1}, and use their algorithm to obtain the same border rank lower bounds as the ones given in Proposition 5 and Theorem 1.
We also compute the $E_{111}$ space for the unique $GL_4$-invariant tensor $T_4\in(\C^4)^*\otimes(\C^4)^*\otimes \Lambda^2\C^4$.
%

\subsection{Acknowledgements}
I would like to thank my advisor J.M. Landsberg for reading numerous drafts and providing many helpful comments.
I thank Giorgio Ottaviani for helpful comments and discussion towards a proof of Theorem 1, and Keller VandeBogert who pointed out the reference \cite{13} for a proof of Theorem 1.
I also thank Laurent Manivel for useful discussions and comments, and Austin Conner for providing an explicit border rank decomposition for the tensor $T_4$ in $\S5.2$.

\section{Preliminaries}

\subsection{Representations of $GL(V)$}

Throughout, let $V$ be a complex vector space of dimension $k$, and let $V^*$ denote the dual space.
The general linear group $GL(V)$ is the group of invertible linear maps $V\to V$.
Irreducible polynomial representations of $GL(V)$ are indexed by partitions $\lambda$ with at most $k$ parts, denoted $S_{\lambda}V$.
Two key representations are those given by $\lambda=(p)$, which is a row of size $p$ and corresponds to the $p$-th symmetric power $\Sym^p(V)$, and $\lambda=(1^p)$, which is a column of size $p$ and corresponds to the $p$-th exterior power $\Lambda^p(V)$.

For two representations $S_\lambda V$ and $S_\mu V$, $S_\lambda V\otimes S_\mu V$ decomposes into a direct sum of irreducible $GL(V)$-modules according to the \emph{Littlewood-Richardson Rule} (see \cite[Chapters~5,8]{9}).
We state two important special cases.
The \emph{Pieri rule} states that
\begin{gather*}
  S_\lambda V\otimes \Sym^pV = \bigoplus_{\nu}S_\nu V\text{,}
\end{gather*}
where the sum is over all partitions $\nu$ obtained by adding $p$ cells to $\lambda$, at most one per column.
The \emph{skew Pieri rule}, states that
\begin{gather*}
  S_\lambda V\otimes \Lambda^pV = \bigoplus_{\nu}S_\nu V\text{,}
\end{gather*}
where the direct sum is over all partitions $\nu$ obtained by adding $p$ cells to $\lambda$, at most one per row.

\subsection{Image of the map $\phi\colon V\to\Hom(S_\mu V,S_\nu V)$}

Following \cite[Section~4]{4}, let $\lambda=(1)$, and $\nu$ be obtained from $\mu$ by adding one cell anywhere not in the first or $k$-th row.
The tensor $T\in V^*\otimes S_\mu V^*\otimes S_\nu V$ corresponds to the $GL(V)$-equivariant map
\begin{gather*}
  \phi\colon V\to\Hom(S_\mu V,S_\nu V)\text{.}
\end{gather*}
Let $\phi_v$ denote the image of $v\in V$ under $\phi$, and $L=\langle v\rangle$ denote the line given by the span of $v$.
Choose a decomposition $V=L\oplus H$, where $H$ is a complementary hyperplane to $L$.
Then $S_\mu V$ and $S_\nu V$ decompose as $(GL(H)\times GL(L))$-modules:
\begin{gather*}
  S_\mu V = \bigoplus_{j\geq0}\bigoplus_{\mu\xrightarrow{j}\pi}S_\pi H\otimes L^j\\
  S_\nu V = \bigoplus_{j\geq0}\bigoplus_{\nu\xrightarrow{j}\pi}S_\pi H\otimes L^j\text{,}
\end{gather*}
where the notation $\mu\xrightarrow{j}\pi$ (resp. $\nu$) means remove $j$ cells from $\mu$ (resp. $\nu$), at most one per column.

Let $y$ denote the cell added from $\mu$ to $\nu$, and $x$ the cell directly above $y$.
We say a partition $\pi$ appears in the domain (resp. codomain) of $\phi_v$ if $S_\pi H\otimes L^j$ is the domain (resp. codomain) for some $j\geq0$.
If $\pi$ appears in both the domain and codomain, then the corresponding $(GL(H)\times GL(L))$-module given by $\pi$ appears in the image; these are the partitions $\pi\subset\nu$ such that $x\in\pi$ and $y\notin \pi$.
If $\pi$ only appears in the domain, then $x\notin \pi$, and the corresponding module is in the kernel of $\phi_v$.
Similarly, if $\nu$ appears only in the codomain, then $y\in \pi$ and the corresponding module is in the cokernel of $\phi_v$.
In particular, $\phi_v$ will always have nonempty kernel and cokernel, so no matrix $\phi_v$ is of full rank.
We obtain a $k$-dimensional space of $m\times n$ matrices where every nonzero matrix has rank equal to the dimension of the image described above.

If $\mu$ differs from $\nu$ in the first row, then every partition $\pi$ that appears in the domain also must appear in the codomain, so every nonzero $\phi_v$ is injective.
If $\mu$ differs from $\nu$ is the $k$-th row, then any partition $\pi$ appearing in the codomain that contains the cell $y$ must have $k$ rows, since $y$ is in the $k$-th row.
However, $H$ is $k-1$ dimensional, so any $GL(H)$-module $S_\pi H$ is 0.
Thus the cokernel of any nonzero $\phi_v$ is trivial, and $\phi_v$ is surjective.

\begin{example}
 Let $\mu=(2,1)$ and $\nu=(2,2)$, with cells $y$ and $x$ labeled.
 \begin{figure}[ht!]
   \begin{center}
     $\mu=$
     \ytableaushort
       {\none x,\none}
       * {2, 1}
    \hspace{53pt}
    $\nu=$
   \ytableaushort
     {\none x,\none y}
     * {2, 2}
   \end{center}
 \end{figure}

 Denoting the partitions corresponding to $GL(V)$- and $GL(H)$-modules with $V$ and $H$ superscripts respectively, for any $v\in V$, we have the decompositions
 \begin{alignat*}{5}
     S_{(2,1)}V
     &= (S_{(2,1)}H) &&  \oplus (S_{(2)}H\otimes L)
     &&  \oplus (S_{(1,1)}H\otimes L) && \oplus (S_{(1)}H\otimes L^2)\\[7pt]
     \ydiagram{2,1}^{\ V}
        &=\ \ydiagram{2,1}^{\ H}
        && \oplus\ \ydiagram{2}^{\ H}
        && \oplus\ \ydiagram{1,1}^{\ H}
        && \oplus\ \ydiagram{1}^{\ H}
 \end{alignat*}
 and
 \begin{alignat*}{4}
     S_{(2,2)}V
     &= (S_{(2,2)}H) &&\oplus (S_{(2,1)}H\otimes L)
     &&\oplus (S_{(2)}H\otimes L^2)\\[7pt]
     \ydiagram{2,2}^{\ V}
     &=\ \ydiagram{2,2}^{\ H}
     &&\oplus\ \ydiagram{2,1}^{\ H}
     &&\oplus\ \ydiagram{2}^{\ H} \text{.}
 \end{alignat*}
 Then,
 \begin{alignat*}{3}
   \ker(\phi_v)
   &= (S_{(1,1)}H\otimes L) &&\oplus (S_{(1)}H\otimes L^2)\\[7pt]
   &=\ \ydiagram{1,1}^{\ H} &&\oplus\ \ydiagram{1}^{\ H} \text{,}\\[7pt]
   \Ima(\phi_v)
   &= (S_{(2,1)}H\otimes L) &&\oplus (S_{(2)}H\otimes L^2)\\[7pt]
   &=\ \ydiagram{2,1}^{\ H} &&\oplus\ \ydiagram{2}^{\ H} \text{,}
 \end{alignat*}
 and
 \begin{alignat*}{3}
  \coker(\phi_v) &= S_{(2,2)}H\\[7pt]
  &=\ \ydiagram{2,2}^{\ H} \text{.}
\end{alignat*}

 %

 When $\dim V=3$, we obtain a three-dimensional space of rank five $6\times8$ matrices.

\end{example}

To explicitly construct a space of matrices, let $\{e_1,\ldots,e_{k}\}$ be a basis for $V$.
Identify the bases of $S_\mu V$ and $S_\nu V$ with semistandard Young tableaux $\SSYT(\mu)$, filled with integers in $\{1,\ldots,k\}$ corresponding to the basis elements of $V$.
Given $\tau\in\SSYT(\mu)$ a basis element of $S_\mu V$, its image under $\phi_{e_i}$, $e_i\in V$, is obtained by adding the cell $y$ to $\tau$, filling $y$ with $i$, acting by the Young symmetrizer, and straightening each tableau in the linear combination (see \cite[Chapters~7,8]{9}).

\begin{example}[$\phi\colon V\to\Hom(S_{(2,1)}V,S_{(2,2)}V)$]\label{ex2}
  Let $V=\C^3$.
 Consider basis elements $e_3\in V$ and $\tau=\
 \begin{ytableau}
   1 & 2\\
   3
 \end{ytableau}\ \in S_{(2,1)}V$.
 Then,
 \begin{align*}
   \phi_{e_3}\left(\ \begin{ytableau}
     1 & 2\\
     3
   \end{ytableau}\ \right)
   &=\ \begin{ytableau}
   1 & 2\\
   3 & 3
 \end{ytableau}\
 + \ \begin{ytableau}
 2 & 1\\
 3 & 3
\end{ytableau}\
 - \ \begin{ytableau}
 3 & 2\\
 1 & 3
\end{ytableau}\
 - \ \begin{ytableau}
 2 & 3\\
 1 & 3
\end{ytableau}\\
&= \ \begin{ytableau}
1 & 2\\
3 & 3
\end{ytableau}\
+ \ \begin{ytableau}
1 & 2\\
3 & 3
\end{ytableau}\
+ \ \begin{ytableau}
1 & 2\\
3 & 3
\end{ytableau}\ - 0\\
&= 3\ \begin{ytableau}
1 & 2\\
3 & 3
\end{ytableau}\ \text{.}
\end{align*}

After normalization, this corresponds to putting $x_3$ in the entry of the matrix indexed by the column of \ \begin{ytableau}
1 & 2\\
3
\end{ytableau}\ and the row of \ \begin{ytableau}
1 & 2\\
3 & 3
\end{ytableau}\ .
This entry is underlined in the matrix below.
Doing this for all $e_i\in V$ and $\tau\in\SSYT(\lambda)$ and normalizing constants, we have the following three-dimensional linear space of rank five $6\times 8$ matrices:
\begin{gather*}
  \begin{pmatrix}
    x_2 & 0   & -x_1 & 0    & 0    & 0    & 0    & 0    \\
    x_3 & x_2 & 0    & -x_1 & -x_1 & 0    & 0    & 0    \\
    0   & x_3 & 0    & 0    & 0    & -x_1 & 0    & 0    \\
    0   & 0   & x_3  & 0    & -x_2 & 0    & -x_1 & 0    \\
    0   & 0   & 0    & \underline{x_3} & 0     & -x_2 & 0    & -x_1 \\
    0   & 0   & 0    & 0    & 0    & 0    & x_3  & -x_2
  \end{pmatrix}\text{.}
\end{gather*}

\end{example}


\subsection{Young flattenings}

Given a tensor $T\in A\otimes B\otimes C$, consider representations $S_{\alpha} A$, $S_{\beta} B$, and $S_{\gamma} C$.
Write the identity map for $S_{\alpha} A$ as $\Id_{S_{\alpha} A}\in S_{\alpha} A^*\otimes S_{\alpha}A$, and similarly for $\Id_{S_{\beta} B}\in S_{\beta} B^*\otimes S_{\beta}B$ and $\Id_{S_{\gamma} C}\in S_{\gamma} C^*\otimes S_{\gamma}C$.
Consider
\begin{gather*}
  T\otimes \Id_{S_{\alpha} A}\otimes \Id_{S_{\beta} B}\otimes \Id_{S_{\gamma} C}
  \in A\otimes B\otimes C\otimes S_{{\alpha}} A^*\otimes S_{{\alpha}}A\otimes S_{\beta} B^*\otimes S_{\beta}B \otimes S_{\gamma} C^*\otimes S_{\gamma}C\text{.}
\end{gather*}
Decompose $S_{\alpha} A\otimes A$ according to the Pieri rule, and project down to an irreducible component $S_{\tilde{\alpha}}A$ given by some partition $\tilde{\alpha}$ obtained by adding a cell to $\alpha$.
Do the same procedure with $S_{\beta} B\otimes B$ and obtain $S_{\tilde{\beta}}B$, and decompose $S_{\gamma} C^*\otimes C$ to obtain $S_{\hat{\gamma}}C^*$, where now $\hat{\gamma}$ is obtained by removing a cell from $\gamma$.
Doing this projection for all three factors gives a tensor
\begin{gather*}
  T' = S_{\alpha}A^*\otimes S_{\beta}B^*\otimes S_{\hat{\gamma}}C^*
  \otimes S_{\tilde{\alpha}} A\otimes S_{\tilde{\beta}} B\otimes S_{\gamma}C\text{.}
\end{gather*}
Consider this tensor as the composition of linear maps
\begin{align*}
  T'\colon S_{\alpha} A\otimes S_{\beta} B\otimes S_{\gamma} C^*
  &\to S_{\alpha} A\otimes S_{\beta} B\otimes S_{\hat{\gamma}}C^*
  \otimes C^*\\[7pt]
  &\to S_{\alpha} A\otimes S_{\beta} B\otimes S_{\hat{\gamma}}C^*
  \otimes A\otimes B\\[7pt]
  &\to S_{\tilde{\alpha}}A\otimes S_{\tilde{\beta}}B\otimes S_{\hat{\gamma}}C^*\text{,}
\end{align*}
called a \emph{Young flattening}.
The border rank lower bound obtained from the Young flattening $T'$ is
\begin{gather*}
  \left\lceil\frac{\rk(T')}{r_1r_2r_3}\right\rceil
\end{gather*}
where $r_1$, $r_2$, and $r_3$ are the generic ranks of matrices in
\begin{gather*}
  A\subset \Hom(S_\alpha A,S_{\tilde{\alpha}}A)\\
  B\subset \Hom(S_\beta B,S_{\tilde{\beta}}B)\\
  C\subset \Hom(S_\gamma C^*,S_{\hat{\gamma}}C^*)
\end{gather*}
respectively.
The special case of ${\alpha}=(1^p)$, ${\beta}=\varnothing$, and ${\gamma}=(1)$, known as \emph{Koszul flattenings}, were used to show border rank lower bounds for various tensors.
See \cite{3},\cite{2} for more discussion.

\section{Border rank lower bounds}

\subsection{Cartan product inclusions}

\Cprod*

For $G=GL(V)$, this proposition gives a nontrivial upper bound for all but one case.
For the $GL_2$-invariant tensor $T\in(\C^2)^*\otimes(\C^2)^*\otimes S_2\C^2$, the maximum possible rank is three, the rank given by the proposition, and this is the only case where the conciseness lower bound matches the upper bound.
In the next smallest example, the proposition says the $GL_2$-invariant tensor $T\in(\C^2)^*\otimes(S_2\C^2)^*\otimes S_3\C^2$ has rank four, but there are rank five tensors in the ambient space (\cite[Proposition~10.3.4.3]{11}).

\begin{proof}
  We flatten $T$ as the space of matrices $\phi\colon C^*\to A^*\otimes B^*$.

  Let $a$ denote a highest weight vector of $A^*$, and similarly $b$ and $c$ for $B^*$ and $C^*$, scaled so that $\phi(c)=a\otimes b$.
  The $G$-orbit of $c$ is a homogeneous variety $X\subset\P C^*$, and any nonzero matrix in $\phi(\hat{X})$ is of rank one.

  Letting $\dim C=n$, we just need to find $n$ points $c_1,\ldots,c_n\in \hat{X}$ such that their span is the entire ambient space $C^*$, which is always possible.
  Then, one rank decomposition is
  \begin{gather*}
    T=\sum_{i=1}^{n}\phi(c_i)\otimes x_i
  \end{gather*}
  where $x_i\in C$ is dual to $c_i\in C^*$.
\end{proof}

The set of such rank $n$ decompositions of $T$ has dimension at least $n\dim X$.

\begin{example}
  We give an explicit decomposition for the $GL_2$-invariant tensor $T\in(S^2\C^2)^*\otimes (S^2\C^2)^*\otimes S^4\C^2$, expressed in bases as the space of matrices
  \begin{gather*}
    \phi((S^4\C^2)^*)=\left\{
    \begin{pmatrix}
      x_0 & x_1 & x_2\\
      x_1 & x_2 & x_3\\
      x_2 & x_3 & x_4
    \end{pmatrix}
    \Bigg|
    x_0,x_1,x_2,x_3,x_4\in\C
    \right\}
    \subseteq (S^2\C^2)^*\otimes (S^2\C^2)^*\text{.}
  \end{gather*}
  The orbit $X$ of the highest weight line is the degree four rational normal curve, which gives the following set of rank one matrices in $\phi((S^4\C^2)^*)$:
  \begin{gather*}
    \phi(\hat{X})=
    \left\{
    \begin{pmatrix}
      s^4    & s^3t   & s^2t^2\\
      s^3t   & s^2t^2 & st^3\\
      s^2t^2 & st^3   & t^4
    \end{pmatrix}
    \Bigg|
    (s,t)\in\C^2\setminus\{(0,0)\}
    \right\}\text{.}
  \end{gather*}
  Since $S^4\C^2$ is five-dimensional, we need to pick five pairwise distinct points of $\P^1$, and consider their images under the degree four Veronese embedding.
  Letting $\zeta$ be a third root of unity, one choice is $[1:0]$, $[1:1]$, $[1:\zeta]$, $[1:\zeta^2]$, and $[0:1]$.
  The matrices those five points correspond to become our $a_i\otimes b_i$ in the rank decomposition.
  Then, expressed in the original basis, we obtain the following rank five decomposition:
  \begin{align*}
    T&=a_1\otimes b_1\otimes (c_1-c_4)\\[7pt]
    &+ \frac{1}{3}(a_1+a_2+a_3)\otimes(b_1+b_2+b_3)
    \otimes(c_2+c_3+c_4)\\[7pt]
    &+ \frac{1}{3}(a_1+\zeta a_2+ \zeta^2 a_3)\otimes(b_1+\zeta b_2+ \zeta^2 b_3)
    \otimes(\zeta^2 c_2+\zeta c_3+c_4)\\[7pt]
    &+ \frac{1}{3}(a_1+\zeta^2 a_2+ \zeta a_3)\otimes(b_1+\zeta^2 b_2+ \zeta b_3)
    \otimes(\zeta c_2+\zeta^2 c_3+c_4)\\[7pt]
    &+ a_3\otimes b_3\otimes (c_5-c_2)\text{.}
  \end{align*}
\end{example}

\subsection{Young flattenings}

For $T\in V^*\otimes S_{\mu}V^*\otimes S_{\nu}V$, using the notation from $\S 2.3$, we consider the Young flattenings with $\beta=\varnothing$ and $\gamma=(1)$.
For partitions $\alpha$ and $\tilde{\alpha}$, where $\tilde{\alpha}$ is obtained from $\alpha$ by adding a cell, we have the associated Young flattening:
\begin{gather*}
  T_{\alpha,\tilde{\alpha}}\colon S_{\alpha} V^*\otimes S_\nu V^*
  \to S_{\alpha} V^*\otimes V^*\otimes S_\mu V^*
  \to S_{\tilde{\alpha}}V^*\otimes S_\mu V^*
\end{gather*}
where $\nu$ \lq\lq gives\rq\rq one of its cells to $\alpha$, to obtain $\mu$ and $\tilde{\alpha}$.
If $a\otimes b\otimes c\in V^*\otimes S_\mu V^*\otimes S_\nu V$ has rank one, then the Young flattening is
\begin{gather*}
  u\otimes v
  \mapsto v(c)P(u\otimes a)\otimes b\text{,}
\end{gather*}
where $u\in S_{\alpha} V^*$, $v\in S_\nu V^*$, and $P\colon S_{\alpha} V^*\otimes V^*\to S_{\tilde{\alpha}}V^*$ is the $GL(V)$-equivariant projection map.
Since $v(c)$ is a complex number, and the image will lie in $S_{\tilde{\alpha}}V^*\otimes\langle b\rangle$, the dimension of the image of the flattening is equal to the dimension of the image of $P(u\otimes a)$.
The vector $a$ is a fixed element of $V^*$ from the tensor, so we are computing the dimension of the image of the map $u\mapsto P(u\otimes a)$ for $u\in S_{\alpha} V^*$.
This map is $\phi_a\in\phi(V^*)\subset\Hom(S_{\alpha} V^*,S_{\tilde{\alpha}}V^*)$, the map given by $a$ in the space of matrices corresponding to the vector space $V^*$ and partitions $\alpha$ and $\tilde{\alpha}$, also given by our construction from Section 2.2.
This gives a border rank lower bound of
\begin{gather*}
  \underline{R}(T)\geq\left\lceil\frac{\rk(T_{\alpha,\tilde{\alpha}})}{\rk(\phi_a)}\right\rceil\text{.}
\end{gather*}

\subsubsection{Koszul flattenings}

Consider the case when $\alpha=(1^p)$ and $\tilde{\alpha}=(1^{p+1})$, the $p$-th Koszul flattening.
The flattening map is
\begin{gather*}
  T^{\wedge p}\colon \Lambda^{p}V^*\otimes S_\nu V^*
  \to \Lambda^{p+1}V^*\otimes S_\mu V^*\text{.}
\end{gather*}

The domain of $T^{\wedge p}$ contains modules given by adding $p$ cells to $\nu$ according to the skew Pieri rule, and the codomain contains modules given by adding $p+1$ cells to $\mu$ also according to the skew Pieri rule.
Let $y$ be in row $j$, and let $\mu_j$ (resp. $\nu_j$) be the number of cells of the $j$-th row of $\mu$ (resp. $\nu$).

By Schur's Lemma, we have the following containments:
\begin{enumerate}[label=(\roman*)]
\item Possible modules in the image are given by the partitions where the $j$-th part is size $\mu_j+1=\nu_j$.
\item The kernel contains at least the modules given by partitions where the $j$-th part is $\mu_j+2=\nu_j+1$.
\item Finally, the cokernel contains at least the modules given by partitions where the $j$-th part is $\mu_j$.
\end{enumerate}

The following theorem states that the three containments described are all equalities.

\pKoszul*

\begin{remark}
  These maps appear as differential maps of \emph{Pieri resolutions} over exterior algebras. These resolutions fit into a broader picture of \emph{Boij-S{\"o}derberg theory}, looking at the semigroup structure of Betti tables of modules over polynomial rings under entrywise addition.
\end{remark}

\subsubsection{Young flattenings beyond Koszul}

We observe that if we choose $\alpha=\mu$ and $\tilde{\alpha}=\nu$, the Young flattening
\begin{gather*}
  T'\colon S_\mu V^*\otimes S_\nu V^*
  \to S_\mu V^*\otimes V^*\otimes S_\mu V^*
  \to S_{\nu}V^*\otimes S_\mu V^*
\end{gather*}
has isomorphic domain and codomain, so the map $T'$ can possibly be an isomorphism.
By construction, we have $\rk(\phi_v)<\min\{m,n\}$ for $\phi_v\in\phi(V)\subset\Hom(S_\mu V,S_\nu V)$, so if the map were of full rank, we have
\begin{gather*}
  \underline{R}(T)\geq\left\lceil\frac{mn}{\rk(\phi_v)}\right\rceil>\max\{m,n\}\text{.}
\end{gather*}
Then this particular Young flattening, specifically adapted to the tensor, would show that no tensor of this form is of minimal border rank.

These Young flattening maps apply to any tensor $T$ that lies in an ambient space $\C^l\otimes \C^m\otimes\C^n$, and they give new equations for border rank lower bounds.
In certain cases where $m\neq n$ they give better border rank lower bounds beyond those given by Koszul flattenings, see Example 4.

\Yflat*
We state partial results towards the conjecture:
\isotypic*
\begin{proof}
  The first inclusion map into the intermediate space is injective, so the full isotypic component maps to the full isotypic component.
  The second map is surjective, so the full isotypic component maps to the full isotypic component.
\end{proof}

For an isotypic component given by $\pi$ in the domain, the isotypic component may appear with higher multiplicity in $S_\mu V^*\otimes V^*\otimes S_\mu V^*$.
We must show when the isotypic component does appear with higher multiplicity, that the submodules given by the inclusion of $S_\mu V^*\otimes S_\nu V^*$ map nontrivially to $S_\nu V^*\otimes S_\mu V^*$.

We verify Conjecture 2 computationally in small cases:
\smalltensors*

\begin{example}\label{ex5}
  Let $\dim V=3$, $\mu=(2,1)$, and $\nu=(2,2)$, so $T\in V^*\otimes S_{(2,1)}V^*\otimes S_{(2,2)}V\cong \C^3\otimes\C^{8}\otimes\C^{6}$.
  This is the space of matrices of constant rank five in Example \ref{ex2}.

  The transpose of the first Koszul flattening is
  \begin{gather*}
    \ytableausetup{boxsize=1em}
    (T^{\wedge1})^t
    \colon \Lambda^2 V\otimes S_{(2,1)}V
    \to V\otimes S_{(2,2)}V\\
    (T^{\wedge1})^t\colon
    \ydiagram{1,1}\
    \otimes\
    \ydiagram{2,1}\
    \to\
    \ydiagram{1}\
    \otimes\
    \ydiagram{2,2}
  \end{gather*}
  where we use the Young diagrams to represent the irreducible modules.
  By the skew Pieri rule, the domain decomposes as
  \begin{gather*}
    \ydiagram{1,1}\
    \otimes\
    \ydiagram{2,1}\
    \cong\
    \ydiagram{3,2}\
    \oplus\
    \ydiagram{3,1,1}\
    \oplus\
    \ydiagram{2,2,1} \text{,}
  \end{gather*}
  and the codomain decomposes as
  \begin{gather*}
    \ydiagram{1}\
    \otimes\
    \ydiagram{2,2}\
    \cong\
    \ydiagram{3,2}\
    \oplus\
    \ydiagram{2,2,1}\ \text{.}
  \end{gather*}
  By Theorem 1, we have the following:
  \begin{gather*}
    \Ima\ (T^{\wedge1})^t=
    \ydiagram{3,2}\
    \oplus\
    \ydiagram{2,2,1}\
  \end{gather*}
  where the second row of the partitions contain $\mu_2+1=\nu_2=2$ cells,
  \begin{gather*}
    \ker\ (T^{\wedge1})^t=
    \ydiagram{3,1,1}
  \end{gather*}
  where the second row of the partition contains $\mu_2=1$ cell,
  and a trivial cokernel.

  In particular, the map is surjective with $\rk(T^{\wedge1})=18$, and
  \begin{gather*}
    \underline{R}(T)\geq\left\lceil\frac{18}{2}\right\rceil=9\text{.}
  \end{gather*}

  The Young flattening $T'$ is
  \begin{gather*}
    T'\colon
    \ydiagram{2,1}\
    \otimes\
    \ydiagram{2,2}\
    \to\
    \ydiagram{2,2}\
    \otimes\
    \ydiagram{2,1}\ \text{.}
  \end{gather*}
  We compute this to be of full rank $6\cdot8=48$, so we obtain a border rank lower bound of
  \begin{gather*}
    \underline{R}(T)\geq\left\lceil\frac{48}{5}\right\rceil=10\text{.}
  \end{gather*}
  The first Koszul flattening for any tensor in $\C^3\otimes\C^6\otimes\C^8$ can only give a border rank lower bound of at most nine, but the more general Young flattening gives new equations that can test for border rank ten.

  In this example, Proposition 2 does not apply. Drawn as diagrams, the domain and intermediate space of the Young flattening $T'$ decompose as follows:
  \begin{gather*}
    \ydiagram{2,1}\
    \otimes\
    \ydiagram{2,2}\
    \cong\
    \ydiagram{4,3}\
    \oplus\
    \ydiagram{4,2,1}\
    \oplus\
    \ydiagram{3,3,1}\
    \oplus\
    \ydiagram{3,2,2}
  \end{gather*}
  \begin{align*}
    \ydiagram{2,1}\
    \otimes\
    \ydiagram{1}\
    \otimes\
    \ydiagram{2,1}\
    &\cong\
    \ydiagram{5,2}\
    \oplus\
    \ydiagram{5,1,1}\
    \oplus\
    \ydiagram{4,3}\ {}^{\oplus2}\\[7pt]
    &\oplus\
    \ydiagram{4,2,1}\ {}^{\oplus4}
    \oplus\
    \ydiagram{3,3,1}\ {}^{\oplus 3}
    \oplus\
    \ydiagram{3,2,2}\ {}^{\oplus3}\ \text{.}
  \end{align*}
  Every module of the domain $S_{(2,1)}V^*\otimes S_{(2,2)}V^*$ appears with higher multiplicity in the intermediate space $S_{(2,1)}V^*\otimes V^*\otimes S_{(2,1)}V^*$.
\end{example}

\tse*

\begin{proof}
  The lower bound from the flattening in Example 4 is ten, and border rank ten fills the ambient space (\cite[Section~6]{10}), so $\underline{R}(T)=10$.
\end{proof}

\section{The case $\dim V=3$}

When $\dim V=3$, if we wish to obtain a space of matrices of bounded rank from a given $\mu$ with at most three rows, we can only add to the second row of $\mu$.
Cutting off any columns of size three, without loss of generality we may start with $\mu$ a two-row partition such that the first row has strictly more cells than the second row.
Let $a$ be the number of columns of size two, and let $b$ be one less than the number of columns of size one ($b$ is the number of columns of size one in $\nu$).
This gives $\mu=(a+b+1,a)$ and $\nu=(a+b+1,a+1)$, and
\begin{gather*}
  m = \dim(S_\mu V) = \frac{(a+b+3)(a+1)(b+2)}{2}\\
  n = \dim(S_\nu V) = \frac{(a+b+3)(a+2)(b+1)}{2}\text{.}
\end{gather*}

For any nonzero $v\in V$,
\begin{gather*}
  \dim(\Ima\phi_v) = \binom{a+b+4}{3} - \binom{a+3}{3} - \binom{b+3}{3}\\
  \dim(\ker\phi_v) = \binom{a+2}{2}\\
  \dim(\coker\phi_v) = \binom{b+2}{2}\text{.}
\end{gather*}

Let $T_{a,b}$ denote the corresponding tensor, and consider the first Koszul flattening.
By Theorem 1, the image contains the modules given by partitions $(a+b+2,a+1)$ and $(a+b+1,a+1,1)$.
The kernel contains the module given by the partition $(a+b+1,a+2)$ if $a>0$, or is trivial if $a=0$.
The cokernel contains the module given by the partition $(a+b+2,a,1)$ if $b>0$, or is trivial if $b=0$.

Thus
\begin{align*}
  \rk(T_{a,b}^{\wedge1})
  &= \dim S_{(a+b+2,a+1)}V + \dim S_{(a+b+1,a+1,1)}V\\[7pt]
  &= \frac{(a+b+4)(a+2)(b+2)}{2} + \frac{(a+b+2)(a+1)(b+1)}{2}\\[7pt]
  &= m + n + a+b+3
  \text{.}
\end{align*}

\begin{prop}
  The first Koszul flattening improves the border rank lower bound of $T_{a,b}$ beyond $\max\{m,n\}$ only if $|a-b|\leq1$.
  \begin{enumerate}
    \item[(i)] If $a=b$, then $m=n$, and $\underline{R}(T_{a,b})\geq m+a+2$.
    \item[(ii)] If $a+1=b$, then $\max\{m,n\}=n$, and $\underline{R}(T_{a,b})\geq n+\left\lceil\frac{a}{2}\right\rceil+1$.
    \item[(iii)] If $a=b+1$, then $\max\{m,n\}=m$, and $\underline{R}(T_{a,b})\geq m+\left\lceil\frac{b}{2}\right\rceil+1$.
  \end{enumerate}
\end{prop}
\begin{proof}
  We consider three cases $a=b$, $a<b$, and $a>b$.

  If $a=b$, then $m=n$, and the improvement is
  \begin{gather*}
    \left\lceil \frac{a+b+3}{2}\right\rceil
    = \left\lceil \frac{2a+3}{2}\right\rceil
    = a+2\text{.}
  \end{gather*}

  When $a<b$, $n$ is larger, and we compute
  \begin{gather*}
    \frac{\rk(T_{a,b}^{\wedge1})}{2} = n + (a-b+2)\ \frac{a+b+3}{4}\text{.}
  \end{gather*}
  Since $\frac{a+b+3}{4}$ is always positive, we look at $a-b+2$.
  To improve the bound, we need $a-b+2>0$.
  By assumption $a<b$, so the only possible solution is $a+1=b$.
  Then $n$ must be even, and the improvement is
  \begin{gather*}
    \left\lceil (a-b+2)\ \frac{a+b+3}{4}\right\rceil
    = \left\lceil \frac{2a+4}{4}\right\rceil
    = \left\lceil \frac{a}{2}\right\rceil + 1\text{.}
  \end{gather*}

Finally, if $a>b$, then $m$ is larger, and we compute
\begin{gather*}
  \frac{\rk(T_{a,b}^{\wedge1})}{2} = l + (b-a+2)\ \frac{a+b+3}{4}\text{.}
\end{gather*}
To improve the bound, we need $b-a+2>0$.
By assumption $a>b$, so we obtain $a=b+1$.
Then $m$ must be even, and the improvement is
\begin{gather*}
  \left\lceil (b-a+2)\ \frac{a+b+3}{4}\right\rceil
  = \left\lceil \frac{2b+4}{4}\right\rceil
  = \left\lceil \frac{b}{2}\right\rceil + 1\text{.}
\end{gather*}
\end{proof}

For Young flattenings using $\alpha=\mu$ and $\tilde{\alpha}=\nu$, the best possible bound we could get is
\begin{gather*}
  \left\lceil\frac{mn}{\rk(\phi_v)}\right\rceil=\left\lceil\frac{(a+b+3)^2(a+2)(b+2)}{2(a+b+4)}\right\rceil\text{.}
\end{gather*}

We also obtain new equations for border rank lower bounds for tensors that lie in the same ambient space.

Computing the Young flattenings for small values of $a$ and $b$, we obtain the following lower bounds.
\smallab*

These bounds are obtained by verifying that each corresponding Young flattening map is full rank.

\section{The case $V\to \Hom(V,\Lambda^2V)$}\label{section5}

\subsection{Koszul flattenings}

For this section, let
\begin{gather*}
  T_{k} = \sum_{1\leq i<j\leq k}\alpha_i\otimes\alpha_j\otimes(e_i\wedge e_j)
  - \alpha_j\otimes\alpha_i\otimes(e_i\wedge e_j)
  \in V^*\otimes V^*\otimes \Lambda^2V
\end{gather*}
 denote the tensor corresponding to $V\to \Hom(V,\Lambda^2V)$ where $\dim V=k$, $\{e_i\}$ is a basis for $V$, and $\{\alpha_i\}$ is the dual basis.

The $p$-th Koszul flattening is the map
\begin{gather*}
  T_k^{\wedge p}\colon
  \Lambda^p V^*\otimes \Lambda^2 V^*
  \to \Lambda^{p+1}V^*\otimes V^*\text{.}
\end{gather*}
The codomain contains only two modules, those given by the partitions $(1^{p+2})$ and $(2,1^p)$.
Both these modules are also in the domain, so by Theorem 1 the map $T_k^{\wedge p}$ is surjective, and its rank is $\dim(\Lambda^{p+1}V^*\otimes V^*)=k\binom{k}{p+1}$.
The border rank lower bound obtained from the flattening is
\begin{gather*}
  \left\lceil\frac{k\binom{k}{p+1}}{\binom{k-1}{p}}\right\rceil
  =\left\lceil\frac{k^2}{p+1}\right\rceil\text{.}
\end{gather*}

Restricting to an $h$-dimensional subspace $H\subset V^*$ does not improve the bounds.
The map $T_k^{\wedge p}$ is surjective, so the restricted flattening map
\begin{gather*}
  \Lambda^pH\otimes \Lambda^2V^*
  \to \Lambda^{p+1}H\otimes V^*
\end{gather*}
is also surjective.
The codomain has dimension $k\binom{h}{p+1}$, which gives a border rank lower bound of
\begin{gather*}
  \left\lceil\frac{k\binom{h}{p+1}}{\binom{h-1}{p}}\right\rceil
  =\left\lceil\frac{kh}{p+1}\right\rceil\text{.}
\end{gather*}
Maximize $\frac{h}{p+1}$, we choose $h=k$ and $p=1$, giving the best lower bound we can obtain from any Koszul flattening:

\begin{restatable}{cor}{skewlb}\label{skewlb}
  $\underline{R}(T_{k})\geq
  \left\lceil\frac{k^2}{2}\right\rceil=\binom{k}{2}+\left\lceil\frac{k}{2}\right\rceil$.
\end{restatable}

\begin{remark}
The $p=1$ Koszul flattening is the same map as the Young flattening constructed for this tensor described in $\S3.2.2$.
\end{remark}

\begin{prop}
  $\underline{R}(T_3) = 5$.
\end{prop}
This was already known, see \cite[Proposition~3.2]{7}.
\begin{proof}
  From Corollary 1,
  \begin{gather*}
    \underline{R}(T_3)\geq5\text{,}
  \end{gather*}
  and border rank five fills the ambient space (\cite[Theorem~4.6]{6}).
\end{proof}

\begin{prop}
  $\underline{R}(T_4) = 8$.
\end{prop}
\begin{proof}
  From Corollary 1,
  \begin{gather*}
    \underline{R}(T_4)\geq8\text{,}
  \end{gather*}
  and border rank eight fills the ambient space (\cite[Corollary~3.10]{6}).
\end{proof}

\subsection{Upper bounds}

We give an explicit border rank five decomposition for $T_3$:
\begin{align*}
  T_3 = \lim_{t\to0}\frac{1}{t}\big[
  &t\alpha_1\otimes \alpha_2\otimes(e_1\wedge e_2)\\
  - &t\alpha_2\otimes\alpha_1\otimes(e_1\wedge e_2)\\
  + &(\alpha_{3}+t\alpha_1)
  \otimes(\alpha_{3}-t\alpha_1)
  \otimes(e_1\wedge e_{3})\\
  + &(\alpha_{3}+t\alpha_2)
  \otimes(\alpha_{3}-t\alpha_2)
  \otimes(e_2\wedge e_{3})\\
  - &\alpha_3\otimes\alpha_3\otimes(e_1\wedge e_3+e_2\wedge e_3)\big]\text{.}
\end{align*}

We thank Austin Conner for the border rank eight decomposition of $T_4$:
\begin{align*}
  T_4 = \lim_{t\to0}
  \frac{1}{4t^3}\Bigg[
  &\left(\alpha_1+2\alpha_3\right)
  \otimes\alpha_1
  \otimes\left(e_2\wedge e_4-2t^3e_1\wedge e_3\right)\tag{1}\\[7pt]
  +&\alpha_1
  \otimes\left(\alpha_1+2\alpha_3\right)
  \otimes\left(e_2\wedge e_4+2t^3e_1\wedge e_3\right)\\[7pt]
  -&4\left(\alpha_1+\alpha_3+t\alpha_2+t^2\alpha_4\right)
  \otimes\left(\alpha_1+\alpha_3-t\alpha_2+t^2\alpha_4\right)
  \otimes\left(e_2\wedge e_4-t^2(e_1\wedge e_2+e_1\wedge e_4)\right)\\[7pt]
  +&2\left(\alpha_1+\alpha_3+2t\alpha_2+2t^2\alpha_4\right)
  \otimes\left(\alpha_1+\alpha_3-2t\alpha_2+2t^2\alpha_4\right)
  \otimes\left(e_2\wedge e_4-t^2(2e_1\wedge e_2+e_1\wedge e_4)\right)\\[7pt]
  -&2\left(\alpha_3-2t\alpha_2\right)
  \otimes\left(\alpha_3+2t\alpha_2\right)
  \otimes\left(e_2\wedge e_4+2t^2(e_1\wedge e_2+e_2\wedge e_3)\right)\\[7pt]
  +&4\left(\alpha_3-t\alpha_2\right)
  \otimes\left(\alpha_3+t\alpha_2\right)
  \otimes\left(e_2\wedge e_4+t^2(e_1\wedge e_2+e_1\wedge e_4+e_2\wedge e_3-e_3\wedge e_4)\right)\\[7pt]
  -&2\left(\alpha_1+\alpha_3+2t\alpha_4\right)
  \otimes\left(\alpha_1+\alpha_3-2t\alpha_4\right)
  \otimes\left(t^2e_1\wedge e_4\right)\\[7pt]
  +&4\left(\alpha_3-t(\alpha_2+\alpha_4)\right)
  \otimes\left(\alpha_3+t(\alpha_2+\alpha_4)\right)
  \otimes\left(t^2(e_3\wedge e_4-e_1\wedge e_4)\right)\Bigg]\text{.}
\end{align*}

We observe that $T_{k}$ splits into two sums:
\begin{gather*}
 T_{k} = \sum_{1=i<j=k-1}\alpha_i\otimes\alpha_j\otimes(e_i\wedge e_j)
 - \alpha_j\otimes\alpha_i\otimes(e_i\wedge e_j)
 + \sum_{j=1}^{k-1}\alpha_j\otimes\alpha_{k}\otimes(e_j\wedge e_{k})
 - \alpha_{k}\otimes\alpha_j\otimes(e_j \wedge e_{k})\\
 = T_{k-1}
 + \sum_{j=1}^{k-1}\alpha_j\otimes\alpha_{k}\otimes(e_j\wedge e_{k})
 - \alpha_{k}\otimes\alpha_j\otimes(e_j \wedge e_{k})
\end{gather*}
where the first sum involves only the first $k-1$ basis elements, and the second one has all terms that use the $k$-th basis element.

Notice that
\begin{gather*}
  T_{k}-T_{k-1}
  =\sum_{j=1}^{k-1}\alpha_j\otimes\alpha_{k}\otimes(e_j\wedge e_{k})
  - \alpha_{k}\otimes\alpha_j\otimes(e_j \wedge e_{k})\\
  = \lim_{t\to 0}\frac{1}{t}\left[
  \sum_{j=1}^{k-1}
  (\alpha_{k}+t\alpha_j)
  \otimes(\alpha_{k}-t\alpha_j)
  \otimes(e_j\wedge e_{k})
  -\alpha_{k}\otimes\alpha_{k}\otimes
  \left(\sum_{i=1}^{k-1} e_i\wedge e_{k}\right)
  \right]
\end{gather*}
which shows $\underline{R}(T_{k}-T_{k-1})=k$.

\begin{theorem}
  For $k\geq4$,
  \begin{gather*}
    \underline{R}(T_{k})\leq\binom{k+1}{2}-2=\binom{k}{2}+k-2\text{.}
  \end{gather*}
\end{theorem}
\begin{proof}
  We induct on $\dim V$.
  Proposition 9 gives the result for $k=4$, our base case.
  Assume $\underline{R}(T_{k-1})\leq\binom{k}{2}-2$.
  We have $T_{k} = T_{k-1} + (T_{k}-T_{k-1})$, so
  \begin{gather*}
    \underline{R}(T_{k})
    \leq\underline{R}(T_{k-1})+\underline{R}(T_k-T_{k-1})\\
    \leq\binom{k}{2}-2+k\\
    =\binom{k+1}{2}-2\text{.}
  \end{gather*}
\end{proof}

\begin{cor}
  \begin{gather*}
    \underline{R}(T_5)=13\text{.}
  \end{gather*}
\end{cor}
\begin{proof}
  An explicit border rank thirteen decomposition of $T_5$ is the border rank eight decomposition (1) of $T_4$ plus the border rank five decomposition of $T_5-T_4$.
  By Corollary 1, thirteen is the lower bound.
\end{proof}

For larger $k$, we have the following bounds:
\begin{gather*}
  18\leq\underline{R}(T_6)\leq19\\
  25\leq\underline{R}(T_7)\leq26\\
  32\leq\underline{R}(T_8)\leq34\\
  41\leq\underline{R}(T_9)\leq43\text{.}
\end{gather*}

\section{Border apolarity}

\subsection{Border apolarity proof of Proposition 5}

\tse*

\begin{proof}[Proof (Border apolarity):]
  We use the $(210)$ border apolarity test given in \cite[Proposition~3.5]{1}.
  We have
  \begin{gather*}
    T_C\colon S_{(2,2)}V^*\to V^*\otimes S_{(2,1)}V^*
    \cong S_{(3,1)}V^*\oplus S_{(2,2)}V^*\oplus S_{(2,1,1)}V^*
  \end{gather*}
  where the domain $S_{(2,2)}V^*$ maps nontrivially to the submodule $S_{(2,2)}V^*$ in the codomain, a six dimensional subspace inside a 24-dimensional space.
  If we wish to prove the border rank is greater than nine, then we must show that for any choice of nine-dimensional Borel-fixed subspace $E_{110}\subset V^*\otimes S_{(2,1)}V^*$ containing $S_{(2,2)}V^*$, the map
  \begin{gather*}
    \psi \colon V^*\otimes E_{110}
    \to \Lambda^2 V^*\otimes S_{(2,1)}V^*
  \end{gather*}
  has a kernel of dimension less than nine.

  Since $E_{110}$ must contain $S_{(2,2)}V^*$, we first determine the rank of the map
  \begin{gather*}
    \psi'\colon V^*\otimes S_{(2,2)}V^*
    \to \Lambda^2 V^*\otimes S_{(2,1)}V^*\text{,}
  \end{gather*}
  the first Koszul flattening.
  In Example \ref{ex5}, the transpose of the first Koszul flattening has trivial cokernel, so $\psi'$ has trivial kernel.

    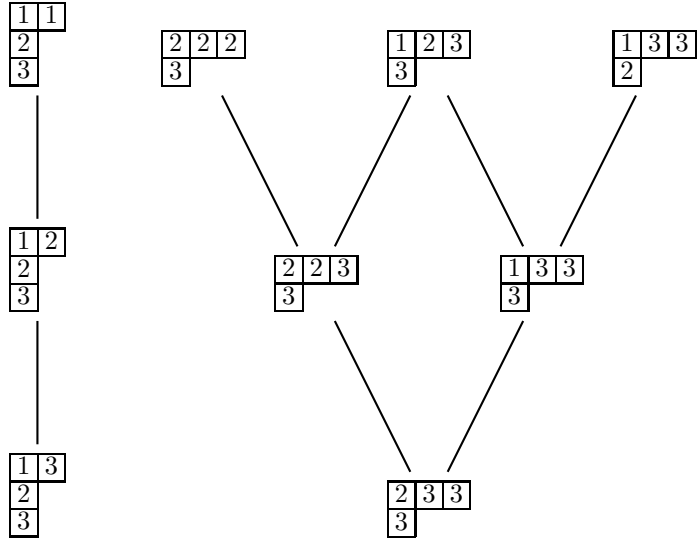
\begin{figure}[ht!]
      \begin{center}
      \begin{tikzpicture}
      \node (top) at (0,6) {\ytableaushort
        {11,2,3}
        * {2, 1, 1}};

      \node (mid) at (0,3) {\ytableaushort
        {12,2,3}
        * {2, 1, 1}};

      \node (bot) at (0,0) {\ytableaushort
        {13,2,3}
        * {2, 1, 1}};

      \draw [thick] (top) -- (mid);
      \draw [thick] (mid) -- (bot);

      \end{tikzpicture}
      \hspace{23pt}
      \begin{tikzpicture}
      \node (top) at (3,-6) {\ytableaushort
        {233,3}
        * {3, 1}};

      \node (mid1) at (1.5,-3) {\ytableaushort
        {223,3}
        * {3, 1}};

      \node (mid2) at (4.5,-3) {\ytableaushort
        {133,3}
        * {3, 1}};

      \node (bot1) at (0,0) {\ytableaushort
        {222,3}
        * {3, 1}};

      \node (bot2) at (3,0) {\ytableaushort
        {123,3}
        * {3, 1}};

      \node (bot3) at (6,0) {\ytableaushort
        {133,2}
        * {3, 1}};

      \draw [thick] (top) -- (mid1);
      \draw [thick] (top) -- (mid2);
      \draw [thick] (mid1) -- (bot1);
      \draw [thick] (mid1) -- (bot2);
      \draw [thick] (mid2) -- (bot2);
      \draw [thick] (mid2) -- (bot3);

      \end{tikzpicture}

      \end{center}
      \caption{The three weight spaces of $S_{(2,1,1)}V^*$ (on the left) and bottom three levels of weight spaces of $S_{(3,1)}V^*$ (on the right).}
    \end{figure}

  Next, we consider the three-dimensional complement $D$ of $S_{(2,2)}V^*$ in $E_{110}$.
  To produce all possible Borel-fixed subspaces, we have four choices of how many basis elements of $S_{(2,1,1)}V^*$ to use.
  If we use all three, then $D=S_{(2,1,1)}V^*$, and the kernel of $\psi$ is three-dimensional.
  If we use the bottom two lowest weight vectors, we can only use the lowest weight vector of $S_{(3,1)}V^*$, and $\psi$ has a four-dimensional kernel.
  If we only use the lowest weight vector of $S_{(2,1,1)}V^*$, there are two choices of two-dimensional Borel-fixed subspaces of $S_{(3,1)}V^*$, shown in Figure 1 by going either left or right from the lowest weight vector.
  Either choice gives a six-dimensional kernel.
  Finally, there are three choices of subspaces if $D\subset S_{(3,1)}V^*$.
  Taking the left three or right three weight vectors gives a six-dimensional kernel, while taking the bottom triangle gives a seven-dimensional kernel.

  Since the kernel of $\psi$ has dimension less than nine for every choice of $E_{110}$, $T$ cannot have border rank nine, so $\underline{R}(T)\geq10$.
  Border rank ten fills the ambient space (\cite[Section~6]{10}), so
  \begin{gather*}
    \underline{R}(T)=10\text{.}
  \end{gather*}
\end{proof}

\subsection{Border apolarity proof of Corollary 1}

\skewlb*

\begin{proof}[Proof (Corollary 1, border apolarity):]
  We again use the $(210)$ border apolarity test.
  Consider a Bored-fixed subspace $E_{110}\subset V^*\otimes V^*$ containing $\Lambda^2 V^*$ of dimension $\binom{k}{2}+d$, and the map
  \begin{gather*}
    \psi
    \colon V^*\otimes E_{110}
    \to \Lambda^2 V^*\otimes V^*\text{.}
  \end{gather*}
  On just $\Lambda^2 V^*\subseteq E_{110}$, we have
  \begin{gather*}
    \psi'
    \colon V^*\otimes \Lambda^2 V^*
    \to \Lambda^2 V^*\otimes V^*\text{,}
  \end{gather*}
  the first Koszul flattening, which we showed to be full rank.

  Thus, for the full $E_{110}$, the kernel of $\psi$ has dimension $kd$.
  For the $(210)$ border apolarity test to fail, we want the kernel to be of smaller dimension than that of $E_{110}$, so solving for $d$ in $\binom{k}{2}+d>kd$ yields
  \begin{gather*}
    d<\frac{k}{2}\leq\left\lceil\frac{k}{2}\right\rceil\text{.}
  \end{gather*}
  Thus,
  \begin{align*}
    \underline{R}(T_k)\geq\binom{k}{2}+\left\lceil\frac{k}{2}\right\rceil\\[7pt]
    = \left\lceil\frac{k^2}{2}\right\rceil\text{.}
  \end{align*}
\end{proof}

\subsection{A Borel-fixed $E_{111}$ space for $T_4$}

We also apply the border apolarity alogrithm for $r=8$ up to multidegree $(111)$ for $T_4$, corresponding to the space of $4\times4$ skew-symmetric matrices.
The eight-dimensional $E_{110}\subset V^*\otimes V^*$ contains the six-dimensional $\Lambda^2 V^*$, so we add two dimensions from $\Sym^2 V^*$.
By the Borel-fixed condition, the only possible addition is $\langle \alpha_4^2,\alpha_3\alpha_4\rangle$.

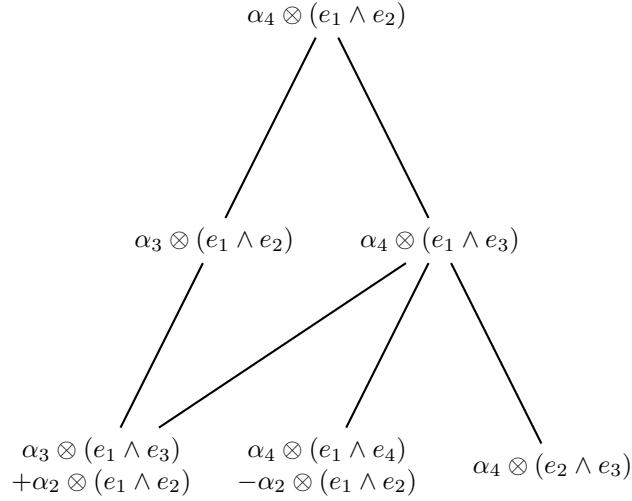
\begin{figure}[ht!]
\begin{center}
\begin{tikzpicture}
\node (top) at (3,6) {$\alpha_4\otimes(e_1\wedge e_2)$};

\node (mid1) at (1.5,3) {$\alpha_3\otimes (e_1\wedge e_2)$};

\node (mid2) at (4.5,3) {$\alpha_4\otimes (e_1\wedge e_3)$};

\node[align=center] (bot1) at (0,0) {$\alpha_3\otimes (e_1\wedge e_3)$\\ $+\alpha_2\otimes(e_1\wedge e_2)$};

\node[align=center] (bot2) at (3,0) {$\alpha_4\otimes (e_1\wedge e_4)$\\ $-\alpha_2\otimes(e_1\wedge e_2)$};

\node (bot3) at (6,0) {$\alpha_4\otimes (e_2\wedge e_3)$};

\draw [thick] (top) -- (mid1);
\draw [thick] (top) -- (mid2);
\draw [thick] (mid1) -- (bot1);
\draw [thick] (mid2) -- (bot1);
\draw [thick] (mid2) -- (bot2);
\draw [thick] (mid2) -- (bot3);

\end{tikzpicture}

\end{center}
\caption{The top three levels of weight spaces of $S_{(1,1,0,-1)}V$.}
\end{figure}

Next, the space $E_{011}\subset V^*\otimes \Lambda^2 V$ must contain $V$, and also has a Borel-fixed four-dimensional subspace of $S_{(1,1,0,-1)}V$, the complement of $V$ in $V^*\otimes \Lambda^2 V$.
There is a one-dimensional family of subspaces that passes both of the $(012)$ and $(021)$ tests:
\begin{align*}
  E_{011}=V\oplus\langle &\alpha_4\otimes e_1\wedge e_2,\\
  &\alpha_3\otimes (e_1\wedge e_2),\\
  &\alpha_4\otimes (e_1\wedge e_3),\\
  &s\alpha_4\otimes (e_1\wedge e_4)
  +s\alpha_3\otimes (e_1\wedge e_3)
  +t\alpha_4\otimes (e_2\wedge e_3)
  \mid[s:t]\in\P^1\rangle\text{.}
\end{align*}
By symmetry, the admissible $E_{101}$ spaces are the same after exchanging the two copies of $V^*$.
Then, we obtain the unique admissible $E_{111}$ space by choosing the point $[1:0]\in\P^1$ for $[s:t]$ parametrizing the one-dimensional family of admissible $E_{011}$ and $E_{101}$ spaces:
\begin{align*}
E_{111}=\langle &T_4,\\
&\alpha_4\otimes\alpha_4\otimes (e_1\wedge e_2),\\
&\alpha_3\otimes\alpha_4\otimes (e_1\wedge e_2),\\
&\alpha_4\otimes\alpha_3\otimes (e_1\wedge e_2),\\
&\alpha_4\otimes\alpha_4\otimes (e_1\wedge e_3),\\
&\alpha_2\otimes\alpha_3\otimes (e_1\wedge e_2)
-\alpha_3\otimes\alpha_2\otimes (e_1\wedge e_2),\\
&\alpha_2\otimes\alpha_4\otimes (e_1\wedge e_2)
-\alpha_4\otimes\alpha_2\otimes (e_1\wedge e_2),\\
&\alpha_3\otimes\alpha_4\otimes (e_1\wedge e_3)
+\alpha_4\otimes\alpha_3\otimes (e_1\wedge e_3)
+\alpha_4\otimes\alpha_4\otimes (e_1\wedge e_4)\rangle\text{.}
\end{align*}

\printbibliography

\end{document}